\documentclass[10pt]{amsart}
\usepackage{amssymb,amsxtra,amsmath,amsfonts}
\usepackage{mathrsfs}
\usepackage{verbatim,enumerate}
\usepackage{graphicx,color,graphics} 
\usepackage{color}
\usepackage[linktocpage=true,colorlinks=true, linkcolor=blue, citecolor=red, urlcolor=blue]{hyperref}

\newtheorem{theorem}{Theorem}
\newtheorem{lemma}[theorem]{Lemma}
\newtheorem{proposition}[theorem]{Proposition}

\theoremstyle{definition}
\newtheorem{definition}[theorem]{Definition}
\newtheorem{example}[theorem]{Example}

\theoremstyle{remark}
\newtheorem{remark}[theorem]{Remark}

\numberwithin{equation}{section}
\numberwithin{theorem}{section}

\newcommand\thref{Theorem \ref}

\newcommand\deref{Definition \ref}

\renewcommand{\comment}[1]{}
\def\CC{\mathbb{C}}

\def\NN{{\mathscr{S}}}
\def\NN{\hat{\mathscr{S}}}

\def\V{\mathcal{V}}

\def\A{\mathcal{A}}

\def\N{\mathscr{S}}

\def\S{\mathscr{S}}

\def\U{\mathcal{U}}
\def\ZZ{\mathbb{Z}}

\def\s{(\mathrm{s})}
\def\n{(\mathrm{n})}

\DeclareMathOperator\End{End}
\DeclareMathOperator\Der{Der}
\DeclareMathOperator\Hom{Hom}

\DeclareMathOperator\ad{ad}

\DeclareMathOperator\gr{gr}

\DeclareMathOperator\Span{span}
\DeclareMathOperator\QF{Fie}

\DeclareMathOperator\F{F}

\def\vac{{\boldsymbol{1}}}  

\def\ii{\mathrm{i}} 

\def\la{\lambda}

\def\ze{\zeta}
\def\z{z}
\def\d{\partial}

\begin{document}

\title{Braided Logarithmic Vertex Algebras}

\author{Bojko Bakalov}
\address{Department of Mathematics,
North Carolina State University,
Raleigh, NC 27695, United States}
\email{bojko\_bakalov@ncsu.edu}

\author{Juan J. Villarreal}
\address{Department of Mathematical Sciences, 
University of Bath, 
Bath BA2 7AY,
United Kingdom}
\email{juanjos3villarreal@gmail.com}


\date{June 12, 2024}



\begin{abstract} We study a family of algebras defined using a locally-finite endomorphism called a braiding map. When the braiding map is semi-simple, the algebra is a generalized vertex algebra, while when the braiding map is locally-nilpotent we have a logarithmic vertex algebra. We describe a method that associates to these algebras non-local Poisson vertex algebras, and we use this relation to build a new example of a generalized vertex algebra motivated by the non-linear Schr\"odinger non-local Poisson vertex algebra.

\end{abstract}

\maketitle


\section{Introduction}\label{sec1}

The following is the main definition of this work (cf.\ \cite{BV}).
\begin{definition}\label{lm}
Let $V$ be a vector space. A \emph{braiding map} $\N$ on $V$ is a linear operator
\begin{equation}\label{logf14}
\N=\sum_{i=1}^{L}\phi_{i}\otimes  \psi_{i}\,, \qquad \phi_{i}, \psi_{i}\in \End(V)  \,,
\end{equation}
satisfying the following two conditions:
\begin{enumerate}

\item\label{2.4-2}
$\N$ is \emph{symmetric}, i.e., $\N=\sum_{i=1}^{L}\phi_{i}\otimes  \psi_{i}=\sum_{i=1}^{L} \psi_{i}\otimes  \phi_{i}
$;

\smallskip
\item\label{2.4-3}

$[\phi_i, \psi_j] = [\phi_i, \phi_j] = [\psi_i, \psi_j] = 0
\,,\qquad 1\le i,j\le L$. 
\end{enumerate}
\end{definition}

Recall that a linear operator $\varphi$ is called \emph{locally finite} if every vector is contained in some finite-dimensional $\varphi$-invariant subspace. A locally finite operator $\varphi$ has a unique Jordan--Chevalley decomposition $\varphi=\varphi^{\s}+\varphi^{\n}$, where $\varphi^{\s}$ and $\varphi^{\n}$ commute, $\varphi^{\s}$ is semi-simple and $\varphi^{\n}$ is locally nilpotent (for every vector $u$ we have $(\varphi^{\n})^{k}u=0$ for some $k\geq 1$). 


Now, we introduce the algebras that we study in this paper, generalizing those of \cite{BV}.

\begin{definition}\label{ba} A \emph{braided logarithmic vertex algebra} (braided logVA)  is a complex vector space $V$ equipped with a vector $\vac\in V$, a braiding map $\N=\sum_{i=1}^{L}\phi_{i}\otimes  \psi_{i}$ where $\phi_i, \psi_i$ are locally finite on $V$, and a family of bilinear products 
\begin{equation*}
\mu_{(n)}\colon V\otimes V\rightarrow V\,, \qquad
a\otimes b \mapsto a_{(n+\S)}b \quad (n\in \ZZ),
\end{equation*}
satisfying the following identities:
\begin{enumerate}
\item[(i)] (\emph{truncation property}) $a_{(n+\S)}b=0$ for $n\gg 0$.
\item[(ii)] (\emph{vacuum axioms}) $\S(\vac\otimes a)=0$, \;  $\vac_{(n+\S)}a=\delta_{n,-1}a$, \;   $a_{(-1+\S)}\vac=a$. 
\item[(iii)]
(\emph{hexagon identity}) 
\begin{equation} \label{eq21}
\N (\mu_{(n)}\otimes I) = (\mu_{(n)}\otimes I) (\N_{13}+\N_{23})\,, \qquad n\in \ZZ\,.
\end{equation}
\item[(iv)]
(\emph{Borcherds identity}) 
\begin{equation}\label{borcherds2}
\begin{split}
\sum_{j\in \ZZ_+} &(-1)^{j}\mu_{(m+n-j)}(I\otimes \mu_{(k+j)})\binom{n+\N_{12}}{j}\\
& -\sum_{j\in \ZZ_+}(-1)^{n+j}\mu_{(n+k-j)}(I\otimes \mu_{(m+j)})(-1)^{\S^{\s}_{12}}\binom{n+\N_{12}}{j} P_{12} \\
& =\sum_{j\in \ZZ_+} \mu_{(m+k-j)}( \mu_{(n+j)}\otimes I) \binom{m+\N_{13}}{j}\,,
\end{split}
\end{equation}
for all $k,m,n\in\ZZ$, where 
$P$ is the transposition on $V\otimes V$ and $\S^{\s}$ is the semi-simple term of the Jordan--Chevalley decomposition of $\S$
(cf.\ \eqref{eq1.3}). 
\end{enumerate}
\end{definition}

The name braided vertex algebras appeared in the literature in quantum vertex algebras; see \cite{EK} (see also \cite[Remark 3.19]{BV}). Note that a braided logVA  with $\S=0$ is the same as a vertex algebra \cite{Bo}. More generally, we have:

\begin{theorem}\label{thm1.4} Let $V$ be a braided logVA with braiding $\S$. Then: 
\begin{enumerate}
\item[(i)] If $\S=\S^{\n}$, then $V$ is a logarithmic vertex algebra as in \cite{BV}.
\item[(ii)] If $\S={\S^{\s}}$, then $V$ is a generalized vertex algebra as in \cite{FFR, DL, BK2}.
\end{enumerate}
\end{theorem}

In \cite{FB}, Poisson vertex algebras were associated to vertex algebras considering (algebraic) limits.  We have the following generalization of this result.

\begin{theorem}\label{thm1.5} Let $V^{\epsilon}$ be a logarithmic vertex algebra over $\CC[\epsilon]/(\epsilon^2)$ with a state-field map $Y^{\epsilon}$ such that $V^0=V^{\epsilon}/\epsilon V_{\epsilon}$ is a commutative vertex algebra. Then $V^0$ naturally acquires the structure of a non-local Poisson vertex algebra. 
\end{theorem}

Recall that a vertex algebra is called \emph{commutative} if all products $\mu_{(n)}=0$ for $n\ge0$; then the product $\mu_{(-1)}$ is commutative and associative \cite{Bo}.
We review the definition of a non-local Poisson vertex algebra (non-local PVA) in Definition \ref{def2.9} below; see \cite{K, DK}. 

In this work, we describe a new example of a generalized vertex algebra in Theorem \ref{thm4.7}, motivated by the non-linear Schr\"odinger non-local PVA from \cite{DK}; see Definition \ref{defscho}. This generalized vertex  algebra is built using an associative algebra with generators $\{ w_{n},  \bar{w}_n \}_{n\in \ZZ}$ subject to the following relations $(m,k\in\ZZ)$:
\begin{equation}\label{eq3.5b}
\begin{split}
\sum_{j\geq 0}w_{m-j}w_{k+j}+w_{k-j}w_{m+j}&=0\, , \\
\sum_{j\geq 0}\bar{w}_{m-j}\bar{w}_{k+j}+\bar{w}_{k-j}\bar{w}_{m+j}&=0 \, , \\
{w}_{m}\bar{w}_{k}-w_{m-1}\bar{w}_{k+1} +\bar{w}_{k}{w}_{m} -\bar{w}_{k-1}{w}_{m+1}&=0\, .
\end{split}
\end{equation}

In \cite{BV2}, we described, in particular,  logarithmic vertex algebras quantizing the potential free boson non-local PVA and the potential Virasoro-Magri non-local PVA. For the non-linear Schr\"odinger non-local PVA the same methods as in \cite{BV2} cannot be used. In this work, the algebraic limits in Theorem \ref{thm1.5} are used to explain the relation between the generalized vertex algebra defined by \eqref{eq3.5b} and the non-linear Schr\"odinger non-local PVA; see \eqref{5.14}. We note that algebraic limits have been also used recently to quantize jet schemes of Poisson vertex algebras in terms of asymptotic algebras of chiral differential operators \cite{AKM,AKM2}.

This paper is organized as follows. In Section \ref{se2}, we use the language of fields to describe braided logVAs; see Theorems \ref{d3.11} and \ref{thm1.4}. In Section \ref{se3}, following \cite{BV},  we prove an existence theorem for braided logVAs.  In Section \ref{se4.1}, we prove Theorem \ref{thm1.5} following ideas of \cite{BV2}. Finally, in Section \ref{se4.2}, we build the generalized vertex algebra defined by relations \eqref{eq3.5b} and we prove that its associated non-local PVA is the non-linear Schr\"odinger PVA.

\section{Logarithmic Fields}\label{se2}



We denote by $V[\![z]\!][ \zeta , z^{-\CC}]$ the space of infinite sums $\sum_{n}f_{n}(\zeta)z^{n}$
 where $f_{n}(\zeta)\in V[\zeta]$ and $n$ runs over the union of finitely many sets of the form $\{-d_{i} + \ZZ_{\geq0}\}$ with $d_i \in  \CC$. 
The space $V [\![z]\!] [\zeta , z^{-\CC}]$ is equipped with the usual action of the derivatives $\partial_{z}$ and $\partial_{\zeta}$, and the operator 
 \begin{equation*}
D_\z := \d_\z+\z^{-1} \d_\ze\,,
\end{equation*}
which acts as the total derivative with respect to $\z$ if we set $\ze=\log\z$. 

Let $\S=\sum_{i=1}^{L}\phi_{i}\otimes  \psi_{i}$ be a braiding map where $\phi_i, \psi_i$ ($1\leq i\leq L$) are locally finite on $V$ with Jordan--Chevalley decompositions 
\[\phi_i=\phi_i^{\s}+\phi_i^{\n},\qquad   \psi_i=\psi_i^{\s}+\psi_i^{\n}\, .\]
Note that, since all the operators $\phi_i, \psi_j$ commute with each other, their semisimple and nilpotent part commute with each other too.
This implies that $\S$ is locally finite on $V\otimes V$ and 
\begin{equation}\label{eq1.3}
\S^{\s}=\sum_{i=1}^{L}\phi^{\s}_{i}\otimes  \psi^{\s}_{i}\,  ,\qquad\S^{\n}=\sum_{i=1}^{L}\phi^{\s}_{i}\otimes  \psi^{\n}_{i}+\phi^{\n}_{i}\otimes  \psi^{\s}_{i}+\phi^{\n}_{i}\otimes  \psi^{\n}_{i}\, .
\end{equation}
We use this decomposition of $\S$ to define the following map on $V\otimes V$:
\[z^{\S}:=z^{\S^{\s}}e^{\zeta \S^{\n}}\, .\]

The next lemma can be easily derived from \eqref{eq1.3} and \deref{lm} \eqref{2.4-3}.

\begin{lemma}\label{lem2.1}
 For any $a, b\in V$, we have 
  \begin{equation}\label{fin}
z^{\S}(a\otimes b)=z^{\S^{\s}}e^{\zeta \S^{\n}}
(a\otimes b)=\sum_{j=1}^{m}\sum_{i=1}^{r} \frac{1}{j!}z^{d_i}\zeta^{j}\bigl(\S^{\n}\bigr)^{j}(a_i\otimes b_i)\, ,
\end{equation}
where $a\otimes b =\sum_{i=1}^r a_i\otimes b_i$ such that $a_i \in  V$ $($respectively, $b_i \in V)$ are common eigenvectors of all $\phi_j^{\s}$ $($respectively, $\psi^{\s}_j)$, $1\leq j\leq L$, and
\begin{enumerate}
\item[(i)] $\S^{\s}(a_i\otimes b_i)=d_i (a_i\otimes b_i)$ \, with \, $d_i\in \CC;$
\item[(ii)] $m\in \ZZ_{\geq 0}$ \, is such that \, $(\S^{\n})^m (a_i\otimes b_i)=0$ \, for all \, $1\leq i\leq r$.
\end{enumerate}
\end{lemma}

\subsection{Quantum fields}

Now we introduce the notion of a quantum field as in \cite{B}.

\begin{definition}\label{d2.1}
Let $V$ be a vector space. The space of \emph{(logarithmic) fields} is defined by
\begin{equation*}
\QF(V) = \Hom(V,V[\![\z]\!][\ze, z^{-\CC}])\, .
\end{equation*}
We denote the elements of $\QF(V)$ as $a(z,\ze)$ or $a(\z)$ for short. 
\end{definition}

On a braided logVA $V$, we define for $a,b\in V$ the series (cf.\ \eqref{fin}) 
\begin{equation}\label{p-modes1}
\begin{split}
Y(a,z)b&=Y(z)(a\otimes b):= \sum_{n\in \ZZ} \mu_{(n)} \bigl( z^{-n-1-\N} (a \otimes b) \bigr) \\
&=\sum_{ \substack{n\in\ZZ } }\sum_{j=1}^{m}\sum_{i=1}^{r} \frac{(-1)^{j}}{j!}  z^{-n-1-d_i}\ze^j \mu_{(n)} \bigl(\bigl(\S^{\n}\bigr)^j (a_i\otimes b_i)\bigr) \,.
\end{split}
\end{equation}
Note that, by the linearity of $\mu_{(n)}$, we have
$a_{(n+\S)}b=\sum_{i=1}^r \mu_{(n)}(a_i\otimes b_i)=\sum_{i=1}^r{a_i}_{(n+\S)}b_i$. 
Now, by Definition \ref{ba} (i)  and \eqref{p-modes1}, the next proposition is obvious.

\begin{proposition}\label{pro2.2} Let $V$ be a braided logVA. Then $Y(a,z)$, defined by \eqref{p-modes1},
is a logarithmic field on $V$ for every $a\in V$.
\end{proposition}

In a braided logVA $V$, we define the \emph{translation operator} $T\in \End(V)$ by 
\begin{equation}\label{tra}
T(a):=a_{(-2+\S)}\vac\,, \qquad a\in V\, .
\end{equation}
Note that $T(\vac)=0$ by definition.

\begin{proposition}\label{pro2.3} For all $a\in V$, we have the \emph{translation covariance} property
\begin{equation}\label{tracov}
[T,Y(a,z)] =D_z Y(a,z) \,. 
\end{equation}
\end{proposition} 
\begin{proof} 
The Borcherds identity \eqref{borcherds2} for $m=0, k=-2$ on $a\otimes b \otimes \vac$ gives:  
\[\mu_{(n)}(a\otimes Tb)-T\mu_{(n)}(a\otimes b)=-n\mu_{(n-1)}(a\otimes b)-\mu_{(n-1)} \bigl(\S( a\otimes b)\bigr)\, .\]
Plugging this into \eqref{p-modes1}, we obtain the translation covariance.
\end{proof}

\begin{proposition}\label{pro2.4} For any braided logVA $V$, the map $Y$ defined in \eqref{p-modes1} satisfies \, $Y(\vac,z)=I$, \, $Y(a,z)\vac\in V[\![z]\!]$, and \, $Y(a,z)\vac\big|_{z=0} = a$. 
\end{proposition}

\begin{proof} 
By the hexagon identity (Definition \ref{ba} (iii)), the components $\phi_i,\psi_i$ of $\S$ are derivations of all products $\mu_{(n)}$. This implies, in particular, that $\phi_i(\vac) = \psi_i(\vac) = 0$, because $\vac_{(-1+\S)}\vac = \vac$ by the vacuum axiom (Definition \ref{ba} (ii)). Then the claims of the proposition follow easily from \eqref{p-modes1} and the vacuum axiom.
\end{proof}

As in \cite{BV}, we will use the notation $z_{ij} = z_i-z_j$, and we introduce the formal variables $\ze_{ij}=\ze_{ji}$, which are thought of as $\log |z_{ij}|$. Additionally, we consider the formal expansions for $n\in\ZZ$ 
\begin{equation}\label{eq2.5}
\begin{split}
\iota_{z_1,z_2} z_{12}^{n+\N} = \iota_{z_1,z_2} z_{12}^{n+{\S^{\s}}} e^{\ze_{12}\N^{\n}}
&:=\sum_{j\in\ZZ_{+}}\binom{n+\N}{j} (-1)^j z_{1}^{n-j+\S^{\s}}z_{2}^{j}\, ,\\
\iota_{z_2,z_1} z_{12}^{n+\N} = \iota_{z_2,z_1} z_{12}^{n+{\S^{\s}}} e^{\ze_{21}{\S^{\n}}}
&:=\sum_{j\in\ZZ_{+}}\binom{n+\N}{j} (-1)^{n+j+{\S^{\s}}}  z_{1}^{j}z_{2}^{n-j+\N}\, .
\end{split}
\end{equation}
In the following proposition, we will utilize the following generalization of delta function introduced in \cite{B}:
\begin{equation}\label{deltaN}
\delta_{\N}(z_1,z_2):= \sum_{n\in\ZZ} z_{1}^{-n-1-\N} z_{2}^{n+\N}\,.
\end{equation}
Note that we have
\begin{equation}\label{deltaNd}
\frac1{j!} D_{z_2}^{j}\delta_{\N}(z_1,z_2)= \sum_{n\in\ZZ}\binom{n+\S}{j} z_{1}^{-n-1-\N} z_{2}^{n-j+\N}\,.
\end{equation}

\begin{proposition} \label{t3.20} Let $V$ be a braided logarithmic vertex algebra. 
For $n\in\ZZ$, we have the following \emph{Borcherds identity}
\begin{equation}\label{borcherds}
\begin{split}
\iota_{z_1,z_2} & Y(z_1)(I\otimes Y(z_2))z_{12}^{n+\S_{12}}  -\iota_{z_2,z_1}Y(z_2)(I\otimes Y(z_1))z_{12}^{n+\S_{12}}   \\
&=\sum_{j\in\ZZ_{+} } \frac1{j!} Y(z_2)(\mu_{(n+j)}\otimes I)D_{z_{2}}^{j}\delta_{\N_{13}}(z_{1},z_{2})\,.
\end{split}
\end{equation}
\end{proposition}

\begin{proof} 
Using the notation $X(z)(a\otimes b):=\sum_{n\in \ZZ} \mu_{(n)}(a \otimes b)z^{-n-1}$ we have that \eqref{borcherds2} is equivalent to 
\begin{equation}\label{borcherdsx}
\begin{split}
& X(z_1)(I\otimes X(z_2)) \sum_{j\in\ZZ_{+}}\binom{n+\N_{12}}{j} (-1)^j z_{1}^{n-j}z_{2}^{j} \\
&-X(z_2)(I\otimes X(z_1))\sum_{j\in\ZZ_{+}}\binom{n+\N_{12}}{j} (-1)^{n+j+\S^{\s}_{12}}  z_{1}^{j}z_{2}^{n-j} P_{12}  \\
&=\sum_{j\in\ZZ_{+} }X(z_2)(\mu_{(n+j)}\otimes I)\sum_{m\in\ZZ}\binom{m+\S_{13}}{j} z_{1}^{-m-1} z_{2}^{m-j}\, .
\end{split}
\end{equation}
We multiply the right-hand side of this identity by $z_1^{-\S_{13}} z_2^{-\S_{23}}$; 
then from \eqref{eq21}, \eqref{eq2.5}, \eqref{deltaNd}, we obtain \eqref{borcherds}.
\end{proof}

Before we state the locality of braided logVAs, we need the following simple lemma. 

\begin{lemma}\label{lem2.7}
 Let $\S$ be as in \eqref{eq1.3}. Then for any $a\in V$, there exist $m\in \ZZ_{\geq 0}$ and an expression
 $a=\sum_{j=1}^{r} a_{j}$ such that 
\begin{enumerate}
\item[(i)] $\phi^{\s}_{i} a_j ={d_{i, j}}a_j$ \, for some \, $d_{i,j}\in\CC$ \, and all \, $1\leq i\leq L$, \, $1\leq j\leq r;$
\item[(ii)] $\bigl(\phi_i^{\n}\bigr)^{m+1} a_j = 0$ \, for all \, $1\leq i\leq L$, \, $1\leq j\leq r$.
\end{enumerate}
Furthermore, for any $b\in V$, we have 
  \begin{equation*}\label{fin2}
z^{\S}(a\otimes b)=\sum_{j=1}^{r}\sum_{k=0}^{mL}\sum_{1\leq i_1, \cdots , i_k\leq L} \frac{1}{k!}\zeta^{k} (\phi^{\n}_{i_1}\cdots \phi^{\n}_{i_k})a_j\otimes z^{\sum _{i=1}^L d_{i,j}\psi_{i}}(\psi_{i_1}\cdots \psi_{i_k}) b\,.
\end{equation*}
\end{lemma}
\begin{proof} 
First, note that all $\phi_i$ commute with each other; hence, all $\phi^{\s}_{i}$ commute and can be diagonalized them simultaneously.
Then the formula for $z^\S$ follows from the identities 
\begin{align*}
z^{\phi^{\s}_{i}\otimes \psi_i}(a\otimes b) &=\sum_{j=1}^r a_j\otimes z^{{d_{i,j}} \psi_i}b\, ,\\ 
z^{\phi^{\n}_{i}\otimes \psi_i}(a\otimes b) &=\sum_{k=0}^m \frac{\zeta^k}{k!}\bigl(\phi^{\n}_{i}\bigr)^k a\otimes (\psi_i)^k b\, .
\end{align*}
\end{proof}

\begin{proposition}\label{cor2.5}  For every $a,b\in V$, there exists $N\in\ZZ_{\geq 0}$ such that for all $c\in V$ 
\begin{equation}\label{eq1}
\begin{split}
Y(z_1)& (I \otimes Y(z_2)) \,\iota_{z_1, z_2} z_{12}^{N+\N} (a\otimes b)\otimes c \\
&= Y(z_2) (I \otimes Y(z_1)) \,\iota_{z_2, z_1} z_{12}^{N+\N} (b\otimes a)\otimes c \,.
\end{split}
\end{equation}
\end{proposition}
\begin{proof}
First, we rewrite \eqref{deltaN} as 
\[\delta_{\N}(z_1,z_2)= \Bigl(\sum_{n\in\ZZ} z_{1}^{-n-1}z_2^{n}\Bigr)(z_1^{-\N} z_{2}^{\N})
= \delta(z_1, z_2) \Bigl(\frac{z_2}{z_1}\Bigr)^{\N} \,,\]
where we use the usual delta function from \cite{K}. Then from the Leibniz rule we have 
\[D_{z_2}^{j}\delta_{\N_{13}}(z_1,z_2)
= \sum_{k=0}^{j} \binom{j}{k} D_{z_2}^{k}\bigl(\delta(z_1, z_2)\bigr)D_{z_2}^{j-k} \Bigl(\frac{z_2}{z_1}\Bigr)^{\N_{13}} \,.\]

From Lemma \ref{lem2.7}, we have the finite sum 
\begin{align*}
&z_1^{-\N} z_{2}^{\N}(a\otimes  c)=\sum_{j=1}^{r}\sum_{l\geq 0}\sum_{1\leq i_1, \cdots , i_l\leq L} \frac{1}{l!}(\zeta_2-\zeta_1)^{l}\\
&\times  (\phi^n_{i_1}\cdots \phi^n_{i_l})(a_j)\otimes z_1^{-\sum _{i=1}^L d'_{i,j}\psi_{i}}z_2^{\sum _{i=1}^L d'_{i,j}\psi_{i}}(\psi_{i_1}\cdots \psi_{i_l}) (c)\, . 
\end{align*}
As we mention before, the finiteness condition of the sum is independent of $c$. 
Hence, there exists $M\in \ZZ_{\geq 0}$ such that if $r\geq M$ then for all $l\in \ZZ_{\geq 0}$ and $c\in V$
\[(\mu_{(r)}\otimes I) (D_{z_2}^{l}(z_1^{-\N_{13}} z_{2}^{\N_{13}})(a\otimes b\otimes c))=0\, .\]
Replacing $N=M$ in identity \eqref{borcherds} gives us the proposition.
\end{proof}

\begin{remark} 
In \cite[Sect.\ 4.3]{BV2}, we found an example where $\S$ is locally finite but can not be expressed as a sum $\sum \phi_i\otimes \psi_i$ where each component $\phi_i, \psi_i$ is locally finite on $V$. In this case, the integer $N$ necessary for the locality of the fields depends on $a,b$ and $c$ (see also \cite[Remark 3.27]{BV}).

\end{remark}

The next theorem gives an equivalent definition of a braided logVA.

\begin{theorem} \label{d3.11}
A braided logarithmic vertex algebra is a vector space $V$, equipped with a vector $\vac\in V$, an endomorphism $T\in \End(V)$,  a linear map 
\[Y\colon  V\to \QF(V)\,, \qquad a\mapsto Y(a,z)\,, \]
and a braiding map $\N$ on $V$,
which are subject to the following axioms:

\medskip
$($i$)$ $Y(\vac,z)=I$, $Y(a,z)\vac\in V[\![z]\!]$, $Y(a,z)\vac\big|_{z=0} = a$,  $T\vac=\S(\vac\otimes a)=0$.

\medskip
$($ii$)$\; 
$[T,Y(a,z)] = \partial_z Y(a,z)$.

\medskip
$($iii$)$\; $\N=\sum_{i=1}^{L}\phi_{i}\otimes  \psi_{i}$ where $\phi_i, \psi_i$ are locally finite on $V$.

\medskip
$($iv$)$\; $\N (Y(z)\otimes I) = (Y(z)\otimes I) (\N_{13}+\N_{23})$.

\medskip
$($v$)$\;  For every $a,b\in V$ there exists $N\in\ZZ_+$ such that $\forall c\in V$ \begin{equation}\label{eq1}
\begin{split}
Y(z_1)& (I \otimes Y(z_2))\iota_{z_1, z_2}(z_{12}^{N+\N}) (a\otimes b)\otimes c \\
&= Y(z_2) (I \otimes Y(z_1))\iota_{z_2, z_1}(z_{12}^{N+\N}) (b\otimes a)\otimes c .
\end{split}
\end{equation}
\end{theorem}

\begin{proof} We have that braided logVAs satisfy the above conditions, by Propositions \ref{pro2.2},  \ref{pro2.3}, \ref{pro2.4} and \ref{cor2.5}. The other direction follows from Proposition \ref{t3.20} in the next section. 
\end{proof}

Now, we prove Theorem \ref{thm1.4}.

\begin{proof}[Proof of Theorem \ref{thm1.4}]

(i) Follows from Theorem \ref{d3.11}, since this is the definition of a logarithmic vertex algebra in \cite{BV}. 

(ii) Following \cite{BK2}, we define $Q, \Delta$ and $\eta$. Let $\{a_i\}_{i\in I}$ be a basis of $V$ such that $\S(a_i\otimes a_j)=s_{ij} a_i\otimes a_j$ for some $s_{ij}\in \CC$. 
First, we define the subsets 
\[{I_i}:=\{j\in I\,   |\,  s_{jl}-s_{il} \in\ZZ \;\; \forall l\in I \}\, .\]
Note that $j\in I_i$ if and only if $I_i =I_j$; hence $I=\bigsqcup I_{i}$. We define 
$Q:=\{I_{i}\}$.
Now, $Q$ has the structure of an abelian group defined by ($r,t\in\ZZ$):
\begin{equation}\label{eq7.3}
r {I_i}+t I_{j}:=\{k\in I\,   |\,  s_{kl}-r s_{il}-t s_{jl} \in\ZZ  \;\; \forall l\in I \} \,.
\end{equation}
For $\alpha=I_{k} \in Q$, we define $V_{\alpha}=\Span\{a_{j} \,|\, j\in\alpha\}$. Then
\[V=\bigoplus_{\alpha\in Q} V_{\alpha}\,. \]

Second, we define $\Delta\colon Q\times Q\rightarrow  \CC/\ZZ$ by
\begin{equation}\label{eq2.4}
\Delta(I_i, I_j):= s_{ij} \mod \ZZ\, .
\end{equation}
Note that the map $\Delta$ is well defined by the construction of $Q$. The map $\Delta$ is symmetric because $\S$ is symmetric, and $\Delta$ is bilinear with respect to the group structure \eqref{eq7.3}. Finally, the map $\eta\colon Q\times Q\rightarrow  \CC^{\times}$ is given by $\eta(I_i, I_j)=1$, the constant map.


It is easy to see that the axioms of a braided logVA give us the axioms of a generalized vertex algebra. 
\end{proof}

\section{Existence of Braided Logarithmic Vertex Algebras}\label{se3}

 This section follows similar steps as in \cite[Section 2]{BV} where the case of locally nilpotent braiding was considered in  detail.
Let ${\S}=\sum_{i=1}^{L}\phi_{i}\otimes  \psi_{i}$ be a braiding map on $V$ as in  Definition \ref{lm}.  Then we have 
the braiding map on $\QF(V)$ given by
\begin{equation}\label{eq3.1}
\NN:=\sum_{i=1}^{L}\ad(\phi_{i})\otimes\ad(\psi_{i}) \, .
\end{equation} 

\begin{definition}\label{def3.1}
We say that a subspace $\V$ of $\QF(V)$ is $\hat{\S}$-\emph{local} if the following two conditions hold: \begin{enumerate}
\item $\ad(\phi_i), \ad(\psi_i)$ are locally finite on $\V$ for $1\le i\le L$.
\item For any $a,b\in \V$ there is some integer $N\ge 0$ such that
\begin{equation}\label{logf11}
\begin{split}
\mu\bigl(\iota_{z_1, z_2}&\bigl(z_{12}^{N+\hat{\S}} \bigr)a(\z_1) \otimes b(\z_2)\bigr) =\mu\bigl( \iota_{z_2, z_1}\bigl(z_{12}^{N+\hat{\S}}\bigr)  b(\z_2)\otimes a(\z_1)\bigr)\,,
\end{split}
\end{equation}
where $\mu$ denotes the composition in $\End(V)$.
\end{enumerate}
\end{definition}

Note that on the $\hat{\S}$-local space $\V$ we can apply Lemma \ref{lem2.1}. 
Next, we define the {$(n+\NN)$-th products} of fields as in \cite{BV}.

\begin{definition}\label{d2.10} Let
$a, b\in \V$ be local logarithmic fields, and $N$ be from \eqref{logf11}.
For $n\in\ZZ$, the \emph{$(n+\NN)$-th product} $\hat{\mu}_{(n)}(a\otimes b)=a_{(n+\NN)}b$ is defined by:
\begin{equation}\label{logf16}
\begin{split}
\bigl(&a_{(n+\NN)}b\bigr)(z,\ze) \\
&:= \frac{D_{\z_1}^{N-1-n}}{(N-1-n)!} \mu \Bigl(\iota_{z_1, z_2}\bigl(z_{12}^{N+\hat{\S}}\, \bigr)  
a(\z_1,\ze_1)\otimes  b(\z_2,\ze_2) \Bigr)\Big|_{ \substack{\z_1=\z_2=\z \\ \ze_1=\ze_2=\ze} } 
\end{split}
\end{equation}
for $n\le N-1$, and $a_{(n+\N)}b:=0$ for $n\ge N$. 
From now on, to simplify the notation, we will suppress the dependence on $\ze$ and understand that setting $\z_1=z$ automatically sets $\ze_1=\ze$.
\end{definition}


It is easy to check that $a_{(n+\NN)}b$ is again a field, and it does not depend on the choice of $N$ satisfying \eqref{logf11}. 
Note that, when $b=I$ is the identity operator, we have $N=0$ and
\begin{equation}\label{logf-aI}
(a_{(-1+\N)}I)(z) = a(z) \,, \qquad (a_{(-2+\N)}I)(z) = D_z a(z) \, .
\end{equation}


Let $\V$ be an $\NN$-local subspace of $\QF(V)$ containing $I$. 
We denote by $\V^{1}$ the minimal subspace of $\QF(V)$ containing $\hat{\mu}_{(n)}(\V\otimes \V)\subset \QF(V)$ for all $n\in \ZZ$, i.e.,
\begin{equation}\label{logf-v1}
\V^{1}=\sum_{n\in \ZZ}\hat{\mu}_{(n)}(\V\otimes \V)\, .
\end{equation}
Note that $D_z\V\subset \V^1$, by \eqref{logf-aI}.
We have the following analog of Dong's Lemma (cf.\ \cite[Lemma 2.14]{BV}).

\begin{lemma}\label{l2.12b}
Let\/ $\NN$ be the braiding map in \eqref{eq3.1} on\/ $\QF(V)$, and\/ $\V$ be an\/ $\NN$-local subspace of\/ $\QF(V)$.
Then the space\/ ${\V^{1}}$ is\/ $\NN$-local.
 \end{lemma}
 
The proof is the same as in \cite[Lemma 2.14]{BV}  and \cite[Corollary 2.21]{BV}, where the case of locally nilpotent braiding was considered.  
 
Applying Lemma \ref{l2.12b} with $\V^1$ in place of $\V$, we construct the local space $\V^2 = (\V^1)^1$. Continuing by induction, we let $\V^{m+1} = (\V^{m})^{1}$, and obtain a sequence of $\NN$-local subspaces
\begin{equation}\label{logf-vbar}
\V=\V^0 \subset \V^{1}\subset \V^2 \subset\cdots ,  \, \qquad \overline{\V}:=\bigcup_{m\geq 0} \V^m\,.
\end{equation}
The union $\overline{\V}$ is $\NN$-local and closed under all $(n+\NN)$-th products.


Now, we assume that our space $V$ is equipped with a vector $\vac\in V$ 
and a linear operator $T\in\End(V)$  
such that $T\vac=0$. 

\begin{definition}\label{dtrcov}
A  field $a\in\QF(V)$ is called \emph{translation covariant} if 
 \begin{equation}\label{logf20}
 [T,a(z)]=D_{z}a(z)\, .
 \end{equation}  
 We denote by $\QF_{T}(V)$ the space of translation covariant fields.  
 \end{definition}

Additionally, we will assume that the braiding map $\NN$ in \eqref{eq3.1} satisfies
  \begin{equation}\label{Tphipsi}
\phi_i \vac=\psi_i\vac=0,\qquad [T,\phi_i]=[T,\psi_i]=0\,, \qquad i=1,\dots,L \,. 
\end{equation}
The following theorem allows us to build examples of braided logVAs.

\begin{theorem}[Existence Theorem]\label{l2.19} 
Let\/ $V$ be a space equipped with a vector\/ $\vac \in V$, a linear operator\/ $T\in \End(V)$ such that\/ $T\vac=0$,
and a braiding map ${\S}=\sum_{i=1}^{L}\phi_{i}\otimes  \psi_{i}$ on\/ $V$ satisfying \eqref{Tphipsi}. 
Set\/ $\NN=\sum_{i=1}^{L}\ad(\phi_{i})\otimes\ad(\psi_{i})$. Suppose that we have an\/ $\NN$-local subspace
of translation covariant fields, 
\begin{equation*}
\V=\Span\bigl\{I,a^{i}(z) \,\big|\, i\in J \bigr\} \subset \QF_{T}(V) \,,
\end{equation*}
where\/ $J$ is an index set,
which is {complete} in the sense that
\begin{equation*}
V=\Span\bigl\{\vac, \, a^{i_{1}}_{(n_1+\N)}\cdots a^{i_{k}}_{(n_k+\N)}\vac \,\big|\, k\ge1, \, i_{1},\dots, i_{k}\in J, \, n_{1},\dots ,n_{k}\in\ZZ \bigr\} \,.
\end{equation*}
Then\/ $\overline{\V} \subset \QF_{T}(V)$, defined by \eqref{logf-vbar}, is an\/ $\NN$-local space containing\/ $\V$ closed under all\/ $(n+\NN)$-th products, and for every\/ $v \in V$ there exists a unique field\/ $Y (v, z)\in \overline{\V}$ such that\/ $Y (v, z)\vac \big|_{z=0}=v$. Moreover, $(V,\vac,T,Y,\N)$ satisfies 
the conditions of
Theorem \ref{d3.11}, i.e., it is a braided logVA.
\end{theorem}
\begin{proof}
The proof follows the same steps in \cite[Theorems 2.29, 3.2]{BV}. 
First, note that the $(n+\NN)$-th product of any two fields $a,b\in\V$ is translation covariant \cite[Lemma 2.25]{BV}:
\begin{equation}\label{eq3.9f}
 a_{(n+\NN)}b \in\QF_{T}(V)\,, \qquad a,b\in\V \,, \; n\in\ZZ
 \end{equation} 
Second, recall that \cite[Lemma 2.26]{BV}:
 \begin{equation}\label{eqpro}
 a(z)\vac\in V[\![z]\!]\,,
 \qquad a\in\QF_{T}(V) \,.
 \end{equation}  
 Hence, the following linear map is well defined:
 \[\Theta\colon \QF_{T}(V)\rightarrow V\, , \qquad \Theta(a):=a(z)\vac|_{z=0}\, .\]


 We define $U:=\Span_{\CC}\{\Theta(a)\,|\,a\in \bar{\V}\}\subset V$, and  the endomorphisms ${a}_{(n+\N)}$ on $U$ by 
\begin{equation}\label{logf24}
\sum_{n\in \ZZ}z^{-n-1}a_{(n+\N)}u:=\mu\Bigl(z^{\sum_{i=1}^{L}\ad(\phi_i)\otimes \psi_i}a(z)\otimes u\Bigr),
\qquad u\in U\,.
 \end{equation}
  \begin{lemma}\label{l2.18} 
  Let\/ $\NN$ be a braiding map as  in \eqref{eq3.1} satisfying \eqref{Tphipsi}.  Then, for any $\NN$-local space $\V\subset \QF_T(V)$, we have 
\begin{equation}\label{Thanb}
\Theta(a_{(n+\NN)}b) = a_{(n+\N)} \Theta(b) \,, \qquad a,b\in\V \,, \;\; n\in\ZZ \,.
\end{equation}
 \end{lemma}
\begin{proof} 
We will prove that 
\begin{equation}\label{aThetab}
\mu\Bigl(z^{\sum_{i=1}^{L}\ad(\phi_i)\otimes \psi_i}a(z)\otimes \Theta(b)\Bigr)  
=\sum_{n\in \ZZ}z^{-n-1}\Theta\bigl(a(\z)_{(n+\NN)} b(\z)\bigr)\,.
\end{equation}
First, note that from \eqref{eqpro}, we have
\begin{equation*}
c(z_1, z_2):=\mu\Bigl(\iota_{z_1, z_2}\bigl(z_{12}^{N+\hat{\S}} \bigr)a(\z_1) \otimes b(\z_2)\Bigr)\vac \in V[\![z_{1},z_{2}]\!][z_{1}^{-\CC}, \ze_{1}]\, .
\end{equation*}
Then, using the locality \eqref{logf11}, we get $c(z_1, z_2)\in V[\![z_{1},z_{2}]\!]$. 

Next, we observe that \[\Theta\bigl(\ad(\psi_i)b(z)\bigr)=[\psi_i, b(z)]\vac\big|_{z=0}=\psi_i \Theta(b)\,.\] 
Hence,
\begin{align*}
\mu\bigl(z_1^{\ad(\phi_i)\otimes \psi_i}a(z_1)\otimes \Theta(b)\bigr)&=\mu\bigl(\iota_{z_1, z_2}(z_{12})^{\NN}a(\z_1) \otimes b(\z_2)\bigr)\vac\big|_{z_2=0}\\
&=\mu\bigl(\iota_{z_1, z_2}(z_{12}^{-N})\bigl(z_{12}^{N+\hat{\S}} \bigr)a(\z_1) \otimes b(\z_2)\bigr)\vac\big|_{z_2=0}\\
&=z_{12}^{-N}c(z_1, z_2)\vac\big|_{z_2=0} \\
&=\sum_{i=1}^{r}z_{1}^{-N}c(z_1, 0)\, ,
\end{align*}
where we used that $\iota_{z_1,z_2} \ze_{12}|_{z_2=0} = \ze_1$ and $\iota_{z_1,z_2} z_{12}^{-N}|_{z_2=0} = z_1^{-N}$. 
Now the coefficient in front of $z_1^{N-n-1}$ in the Taylor expansion of $c(z_1 , 0)$  is
\begin{align*}
\frac1{(N-n-1)!} & D_{z_1}^{N-n-1} c(z_1,0) \big|_{z_1=0} \\
&= \frac1{(N-n-1)!} \bigl( D_{z_1}^{N-n-1} c(z_1,z_2) \big|_{z_1=z_2=z} \bigr) \big|_{z=0} \\
&= \Theta\bigl( {a}_{(n+\NN)}b \bigr),
\end{align*}
where we used \eqref{eq3.9f}. This proves \eqref{aThetab}, which together with \eqref{logf24} for $u=\Theta(b)$, implies \eqref{Thanb}.
\end{proof}

The rest of the proof of Theorem \ref{l2.19} is the same as in \cite[Theorems 2.29, 3.2]{BV} and is omitted. 
\end{proof}


Finally, we state the following two consequence of Theorem \ref{l2.19} and Lemma \ref{l2.18}; cf.\ \cite[Proposition 3.5 and Theorem 3.28]{BV}.

\begin{proposition}[$(n+\N)$-th product identity]\label{pro3.12} 
If $(V,\vac,T,Y,\N)$ satisfies the assumptions of \thref{d3.11}, then  for any $a, b\in V$ we have  
\begin{equation}\label{nprod-1}
Y(a,z)_{(n+\NN)}Y(b,z)=Y(a_{(n+\N)}b,z)\,,
\qquad 
n\in\ZZ \,.
\end{equation}
In particular,
$Y(Ta,\z)=D_{z}Y(a,z)$.
\end{proposition}


\begin{proposition}[Borcherds identity] \label{t3.20} 
If\/ $(V,\vac,T,Y,\N)$ satisfies the assumptions of \thref{d3.11}, then 
the Borcherds identity \eqref{borcherds} holds.
\end{proposition}

\section{Poisson Limits and Non-linear Schr\"odinger Non-local PVA}\label{se4}

In this section, first we prove Theorem \ref{thm1.5}. Then, we build the
generalized vertex algebra defined by \eqref{eq3.5b} and we prove that its associated non-local PVA is
the non-linear Schr\"odinger PVA.


\subsection{Poisson limits}\label{se4.1}
We start by reviewing some definitions; for more details, we refer to \cite{K,DK,BV2}.

\begin{definition}\label{def2.2}
A \emph{non-local Lie conformal algebra} (abbreviated non-local LCA)  is a $\mathbb{C}[\partial]$-module $\mathcal{V}$ equipped with an admissible $\la$-bracket
$$\{ \cdot_{\lambda}\cdot \} \colon\mathcal{V}\otimes\mathcal{V}\rightarrow\mathcal{V}(\!(\lambda^{-1})\!)\, ,$$
 satisfying the following axioms:
\begin{align*}
\text{(sesqui-linearity)}&\quad \{\partial a_{\lambda}b\}=-\lambda \{a_{\lambda}b\}\, ,\quad\{a_{\lambda}\partial b\}=(\lambda+\partial)\{a_{\lambda} b\}\, ,\\
\text{(skew-symmetry)}&\quad \{b_{\lambda}a\}=-\{a_{-\lambda-\partial}b\}\, ,\\
\text{(Jacobi identity)}&\quad \{a_{\lambda}\{b_{\mu}c\}\}=\{\{a_{\lambda}b\}_{\lambda+\mu}c\}+\{b_{\mu}\{a_{\lambda}c\}\}\, .
\end{align*}
The admissibility of a $\la$-bracket is a technical condition, which ensures that the compositions in the Jacobi identity are well defined.
\end{definition}

\begin{definition}\label{def2.9}
A \emph{non-local Poisson vertex algebra} (abbreviated non-local PVA)  is a non-local LCA $\mathcal{V}$, which in addition has a commutative associative product
$ab\in\mathcal{V}$ for $a,b\in\mathcal{V}$, with a unit $1\in\mathcal{V}$, such that $\partial$ is a derivation of the product and the Leibniz rule holds:
\begin{align*}
\text{(derivation)}&\quad \partial(ab) = (\partial a)b + a(\partial b) \,, \\
\text{(Leibniz rule)}&\quad \{a_{\lambda}bc\}=\{a_{\lambda}b\}c+b\{a_{\lambda}c\} \,.
\end{align*}
\end{definition}

 
Now we recall the main result of \cite{BV2}.
 In \cite{BV2}, we consider an increasing \emph{filtration} on a logarithmic vertex algebra $V$ by subspaces
$\{0\}\subset \F^{0}V\subset \F^{1}V\subset \F^{2}V\subset\cdots$ satisfying the following properties: $\vac\in \F^{0}V$, $T(\F^{n}V)\subset \F^{n}V$ and 
\begin{equation}\label{eq4.1}
\begin{split}
&\N(\F^{m}V\otimes \F^{n}V)\subset \F^{m+n-1}(V\otimes V) := \displaystyle\sum_{k=0}^{m+n-1}\F^{k}V\otimes \F^{m+n-1-k}V \,,\\
&\mu_{(j)}(\F^{m}V\otimes \F^{n}V)\subset 
                \begin{cases}
                  \F^{m+n}V\, ,  & j<0 \,, \\
                  \F^{m+n-1}V\, ,  &  j\geq 0\,,
                \end{cases}
  \end{split}
  \end{equation}
for all $m,n\in\ZZ_+$, where $\F^{n}V := \{0\}$ for $n<0$.
Then, by \cite[Theorem 3.3]{BV2}, the \emph{associated graded} $\gr V:=\bigoplus_{n\in\ZZ_+} \gr^{n} V$, where $\gr^{n} V:=\F^{n}V/\F^{n-1}V$,
has the structure of a non-local PVA defined as follows.
The unit $1\in\gr^0 V = \F^{0}V$ is $\vac$; the derivation $\partial$ on $\gr V$ is induced from $T$; the commutative associative product on $\gr V$ is induced from $\mu_{(-1)}$;
and the $\lambda$-bracket on $\gr V$ is given by
%
%
\begin{equation}\label{eq2.1}
\{a_{\lambda}b\}:=\sum_{n\in\ZZ_+}\frac{\lambda^{n}}{n!}\mu_{(n)}(a\otimes b)
+\mu_{(-1)}\Bigl(\N\, \Bigl(\frac{1}{\lambda+\partial }\,  a\otimes b\Bigr)\Bigr). 
\end{equation}
%

In \cite{FB}, Poisson vertex algebras were associated to vertex algebras considering (algebraic) limits. In this section, we follow the same idea. First, we define ideals. An ideal $W$ of a braided logVA $V$ is a subspace invariant under $T,\phi_i, \psi_i$ ($1\leq i\leq L$) such that $\mu_{(n)}(a\otimes b)\in W$ for all $a\in V$, $b\in W$, $n\in \ZZ$. As for ordinary vertex algebras, it is not hard to see that $V/W$ has the structure of a braided logVA.

Next, we note that it is straightforward to define the notion of a braided logVA $V$ over a commutative ring $R$ instead of $\CC$, so that
$V$ is an $R$-module and $T,\psi_i, \phi_i, \mu_{(n)}$ ($1\leq i\leq L$, $n\in \ZZ$) commute with the action of $R$ on $V$. In particular, we will be
interested in the case $R=\CC[\epsilon]/(\epsilon^2)$. 

\begin{example}\label{ex4.1}
If $V$ is a braided logVA, then $V^{\epsilon}:=V\otimes \CC[\epsilon]/(\epsilon^2)$ is a braided logVA over $\CC[\epsilon]/(\epsilon^2)$, the subspace $\epsilon V^{\epsilon}$ is an ideal of $V^{\epsilon}$,
and the quotient $V^{\epsilon} / \epsilon V^{\epsilon} \cong V$.
\end{example}

Now we prove Theorem \ref{thm1.5}.

\begin{proof}[Proof of Theorem \ref{thm1.5}]
Since $Y^0=Y^{\epsilon}$ mod $\epsilon$ is commutative, we have 
\begin{equation}\label{eq4.2}
\begin{split}
\S(V^{\epsilon}\otimes V^{\epsilon})&\subset \epsilon (V^{\epsilon}\otimes V^{\epsilon})\,, \\
\mu_{(j)}(V^{\epsilon}\otimes V^{\epsilon})&\subset 
\epsilon V^{\epsilon}\,, \quad j\geq 0\, .
\end{split}
\end{equation}
For $a,b\in V^{\epsilon}$, we define a $\lambda$-bracket by
\begin{equation*}\label{eq2.1b}
\{a_{\lambda}b\}:=\sum_{n\in\ZZ_+}\frac{1}{\epsilon}\frac{\lambda^{n}}{n!}\mu_{(n)}(a\otimes b)
+\frac{1}{\epsilon}\mu_{(-1)}\Bigl(\N\, \Bigl(\frac{1}{\lambda+\partial }\,  a\otimes b\Bigr)\Bigr)\,\,\, \mod\epsilon\in V^0 (\!(\lambda^{-1})\!) \,.
\end{equation*}
When $a\in\epsilon V^{\epsilon}$ or $b\in\epsilon V^{\epsilon}$, we have $\{a_{\lambda}b\} = 0$; hence, this descends to a well-defined $\lambda$-bracket on the quotient 
$V^0$.

Now we check that the above bracket satisfies each of the conditions of Definitions \ref{def2.2} and \ref{def2.9}. The main point is that the relations \eqref{eq4.2} are similar to \eqref{eq4.1}. Hence, the proof that $V^0$ is a non-local Poisson vertex algebra is similar to the proof in \cite[Theorem 3.3]{BV2}.

For instance, let us check sesqui-linearity. 
From Proposition \ref{pro3.12}, $Y(Ta,z)=D_{z}Y(a,z)$ on $V^\epsilon$; hence $\mu_{(n)}(Ta\otimes b)=-\mu_{(n-1)}(n+\S)(a\otimes b )$. Then from \eqref{eq4.2} on $V^0$, we get that
\begin{align*}
&\frac{1}{\epsilon}\mu_{(n)}(\partial a \otimes b)=-\frac{1}{\epsilon}n\mu_{(n-1)}(a \otimes b) \text{ mod } \epsilon \; \text{ for } \; n\geq 1\, , \\
&\frac{1}{\epsilon}\mu_{(0)}(\partial a \otimes b)=-\frac{1}{\epsilon}n\mu_{(-1)}(\S ( a \otimes b)) \text{ mod } \epsilon \; .
\end{align*}
Then as in the proof of \cite[Proposition 3.5 (i)]{BV2}, we obtain $\{\partial a_{\lambda}b\}=-\lambda\{a_{\lambda}b\}$. 

Similarly, the skew-symmetry 
is proved as in the proof of \cite[Proposition 3.5 (iii)]{BV2}. The Leibniz rule is proved exactly as in \cite[Proposition 3.6]{BV2}, and by \cite[Proposition 3.5]{BV2} the bracket is admissible. Finally, the Jacobi identity is proved as in \cite[Proposition 3.8]{BV2}. Note that identity \cite[(3.17)]{BV2} is satisfied mod $\epsilon$ multiplying by $\frac{1}{\epsilon^2}$.
%
\end{proof}

\subsection{The non-linear Schr\"odinger generalized vertex algebra}\label{se4.2}
We start by reviewing the following example of a non-local PVA from \cite{DK}.

\begin{definition}\label{defscho} The \emph{non-linear Schr\"odinger non-local PVA}
is defined as the algebra of differential polynomials
$\V:=\CC[u,v,\partial u, \partial v, \dots ]$
 with  $\lambda$-bracket given by:  
\begin{equation}\label{eqs1}
\begin{aligned}
\{u_{\lambda}u\}&=v(\lambda+\partial)^{-1}v, & \{u_{\lambda}v\}&=-v(\lambda+\partial)^{-1}u,\\
\{v_{\lambda}u\}&=-v(\lambda+\partial)^{-1}u, & \{v_{\lambda}v\}&=u(\lambda+\partial)^{-1}u\, .
\end{aligned}
\end{equation}
\end{definition}

Note that we can express equivalently the non-linear Schr\"odinger non-local PVA using the basis $w:=u+\ii v$, $\bar{w}:=u-\ii v$. Then $\V=\CC[w, \bar{w}, w', \bar{w}', \dots]$
(where $w':=\partial w$), and the $\lambda$-bracket is given by 
\eqref{eq2.1} for $a,b\in\{w, \bar{w}\}$,
where 
\begin{equation*}
\S={\S^{\s}}:=\Phi\otimes \Phi \,, \quad
\Phi(w)=\ii w \,, \quad \Phi(\bar{w})=-\ii\bar{w} \,,
\end{equation*}
$\mu_{(-1)}$ denotes the product in $\V$, and
\begin{equation}\label{eq4.6}
\mu_{(n)}(a\otimes b)=0,\qquad a,b\in \{w,\bar{w}\}\,, \;\; n\geq 0\,.
\end{equation}

We use the above relations to construct a braided logVA $V$. We will assume that it exists and will derive the properties that $V$ should satisfy.

 
 
\begin{lemma}\label{lem4.5}
Let $V$ be a braided logVA with elements $w,\bar{w}\in V$, whose\/ $(n+\N)$-th products satisfy \eqref{eq4.6}. Suppose that there is\/ $\Phi\in\Der(V)$  acting as $\Phi(w)=\ii w$, $\Phi(\bar{w})=-\ii\bar{w}$, such that $\N=\Phi\otimes \Phi$. Then 
\[Y(w_{}, z)=\sum_{n\in \ZZ}w_{(n+\S)}z^{-n-1-\ii\Phi}\, ,  \qquad Y(\bar{w}_{}, z)=\sum_{n\in \ZZ}\bar{w}_{(n+\S)}z^{-n-1+\ii\Phi}\,  \]
where the linear operators $w_n:=w_{(n+\S)}$, $\bar{w}_n:=\bar{w}_{(n+\S)}$ satisfy relations \eqref{eq3.5b}. 
\end{lemma}
\begin{proof} 
The form of $Y(w_{}, z)$ and $Y(\bar{w}_{}, z)$ follows from \eqref{p-modes1}, using that $\S^{\n}=0$.
From the Borcherds identity \eqref{borcherds2} with $n=0$ on $w\otimes w\otimes v$, $v\in V$, we have:  
\begin{align*}
&\sum_{j\in \ZZ_+}(-1)^{j}\mu_{(m-j)}(I\otimes \mu_{(k+j)})\binom{-1}{j}(w\otimes w\otimes v)=\sum_{j\geq 0}w_{(m-j+\S)}w_{(k+j+\S)}v\, , \\
& \sum_{j\in \ZZ_+}(-1)^{j}\mu_{(k-j)}(I\otimes \mu_{(m+j)})(-1)^{-1}\binom{-1}{j}(w\otimes w\otimes v)=-\sum_{j\geq 0}w_{(m-j+\S)}w_{(k+j+\S)}v\, , \\
&\sum_{j\in \ZZ_+} \mu_{(m+k-j)}( \mu_{(j)}\otimes I) \binom{m+i I\otimes I\otimes \Phi}{j}(w\otimes w\otimes v)=0\, , 
\end{align*}
where in the last identity we used \eqref{eq4.6}. Hence, from \eqref{borcherds2} we obtain the first equation in \eqref{eq3.5b}. 
The other two relations in \eqref{eq3.5b} are derived analogously from \eqref{borcherds2} for $\bar{w}\otimes \bar{w}\otimes v$ and $w\otimes \bar{w}\otimes v$,
respectively.
%
\end{proof}
We now use the previous results to construct $V$. We proceed similarly to the construction in \cite[Section 4.4]{BV}, and \cite[Section 4.3]{BV2}. 
We define a unital associative topological algebra $\A$ with generators $\{ w_{n},  \bar{w}_n \,|\, n\in \ZZ \}$ subject to the relations \eqref{eq3.5b}. 
The algebra $\A$ admits a module $V$, generated by an element $\vac\in V$ such that $w_{n}\vac=\bar{w}_n\vac=0$ for $ n\geq 0$, and $V$ is linearly spanned by monomials of the form 
\begin{equation}\label{eq5.10}
a^{1}_{n_1} \cdots a^{r}_{n_r} \vac \quad\text{where}\quad
n_i \le -1 \,, \; a^i\in\{w,\bar{w}\} \,, \; 1\le i\le r \,.
\end{equation}


Next, we define linear operators $T$ and $\Phi$ on $V$ by $T(\vac)=\Phi(\vac)=0$ and
\begin{equation}\label{eq5.11}
\begin{split}
&[\Phi, w_n]=\ii w_n, \quad [\Phi, \bar{w}_n]=-\ii\bar{w}_n ,\\
&[T, w_{n}]=-nw_{n-1}-\ii w_{n-1}\Phi, \quad [T, \bar{w}_{n}]=-n\bar{w}_{n-1}+\ii\bar{w}_{n-1}\Phi ,
\end{split}
\end{equation}
motivated by the hexagon identity \eqref{eq21} and translation covariance \eqref{tracov}.
Additionally, we define the fields
\begin{equation}\label{wfie}
w_{} (z):=\sum_{n\in \ZZ}w_{n}z^{-n-1-\ii\Phi}\, ,  \qquad \bar{w}_{}( z):=\sum_{n\in \ZZ}\bar{w}_{n}z^{-n-1+\ii\Phi}\, . 
\end{equation}
We let the braiding map $\N=\Phi\otimes \Phi$. As before, this gives rise to a braiding map $\NN:=\ad(\Phi) \otimes\ad(\Phi)$ on the space of fields $\QF(V)$.

\begin{theorem}\label{thm4.7}
The fields\/ $w(z)$, $\bar{w}(z)$ are translation covariant and\/ $\NN$-local, and they generate
the structure of a braided logVA on\/ $V$. 
\end{theorem}
\begin{proof} We verify each of the hypotheses of Theorem \ref{l2.19}. By definition we have $\vac\in V$ such that $T\vac=0$ and $\S(\vac\otimes a)=0$. By the Jacobi identity, $[[T,\Phi],w_n]=[T,[\Phi,w_n]]-[\Phi,[T,w_n]]=0$, and similarly, $[[T,\Phi],\bar{w}_n]=0$. Hence $[T, \Phi]=0$. 
Next, we observe that by construction the fields $w(z), \bar{w}(z)$ are translation covariant (cf.\ \eqref{tracov}).
We will prove that the space
\[\U=\Span\{I, w(z), \bar{w}(z)\}\subset \QF_T(V) \]
is $\NN$-local.
We have $z_1^{-\ii\Phi}w(z_2)=w(z_2)z_1^{1-\ii\Phi}$ and $z_1^{\ii\Phi}w(z_2)=w(z_2)z_1^{-1+\ii\Phi}$.
Then using \eqref{eq3.5b}, we obtain:
\begin{align*}
\mu\bigl(\iota_{z_1, z_2}\bigl(z_{12}^{\hat{\S}} \bigr)w(\z_1) \otimes w(\z_2)\bigr)&=\iota_{z_1, z_2}\bigl(z_{12}^{-1} \bigr)w(\z_1) w(\z_2)\\
&=\sum_{m,k\in \ZZ}\Bigl(\sum_{j\geq 0}w_{m-j} w_{k+j}\Bigr)z^{-m-1-\ii\Phi}_{1}z_2^{-k-1-\ii\Phi}\\
&=- \sum_{m,k\in \ZZ}\Bigl(\sum_{j\geq 0}w_{k-j} w_{m+j}\Bigr)z^{-m-1-\ii\Phi}_{1}z_2^{-k-1-\ii\Phi}\\
&=\iota_{z_2, z_1}\bigl(z_{12}^{-1} \bigr)w(\z_2) w(\z_1)\\
&=\mu\bigl(\iota_{z_2, z_1}\bigl(z_{12}^{\hat{\S}} \bigr)w(\z_2) \otimes w(\z_1)\bigr)\, .
\end{align*}
Hence, we have the locality of $w(z)$ with itself as in \eqref{logf11} for $N=0$. The locality of $\bar{w}(z), \bar{w}(z)$ and $w(z), \bar{w}(z)$  is proved in a similar way.  
\end{proof}


\begin{remark}\label{rem:pbw}
One can show that the braided logVA $V$ in Theorem \ref{thm4.7} has a PBW basis consisting of lexicographically 
ordered monomials of the form \eqref{eq5.10}. The proof is similar to the usual Poincar\'e--Birkhoff--Witt Theorem for Lie algebras
(see, e.g., \cite{J}). In the context of logarithmic vertex algebras, a version of the PBW Theorem was proved in \cite{H}, 
while the case of quadratic associative algebras is discussed in \cite{PP}.
\end{remark}

Finally, following \cite[Section 16.3]{FB}, we deform the braided logVA in Theorem \ref{thm4.7} by introducing a formal parameter $\epsilon$. 
We consider the unital associative topological algebra $\A^{\epsilon}$ with generators $\{ w_{n},  \bar{w}_n \,|\, n\in \ZZ \}$ subject to the relations:  
\begin{equation}\label{5.14}
\begin{split}
\sum_{j\geq 0}(-1)^j\binom{-\epsilon}{j}(w_{m-j}w_{k+j}-w_{k-j}w_{m+j})&=0\, , \\
\sum_{j\geq 0}(-1)^j\binom{-\epsilon}{j}(\bar{w}_{m-j}\bar{w}_{k+j}-\bar{w}_{k-j}\bar{w}_{m+j})&=0 \, , \\
{w}_{m}\bar{w}_{k}-\epsilon w_{m-1}\bar{w}_{k+1} -\bar{w}_{k}{w}_{m} +\epsilon\bar{w}_{k-1}{w}_{m+1}&=0\, .
\end{split}
\end{equation}
As in Theorem \ref{thm4.7}, $\A^{\epsilon}$ has a module $V^{\epsilon}$ spanned by the monomials \eqref{eq5.10}.
We have operators $T, \Phi\in \End_{\CC[\epsilon]}(V^{\epsilon})$,
fields 
\[w^{\epsilon}(z)=\sum w_{n}z^{-n-1}e^{i\epsilon \Phi \zeta}\, , \qquad  \bar{w}^{\epsilon}(z)=\sum \bar{w}_{n}z^{-n-1}e^{-i\epsilon \Phi \zeta}\, ,\]
a braiding map $\N=\epsilon \Phi\otimes \Phi$, and $\NN:=\epsilon\ad(\Phi) \otimes\ad(\Phi)$ on $\QF(V^{\epsilon})$. The same steps as in the proof of Theorem \ref{thm4.7} give us that $V^{\epsilon}$ is a braided logVA. Moreover, $V^0:=V^{\epsilon}/\epsilon V^{\epsilon}$ is commutative. Then Theorem \ref{thm1.5} implies that $V^0$ has the structure of a non-local PVA. It is not hard to see that this non-local PVA is exactly the non-linear Schr\"odinger non-local PVA from Definition \ref{defscho}.


\subsection*{Acknowledgments}
The first author was supported in part by a Simons Foundation grant 584741. 
He is grateful to the University of Bath for the hospitality during August 2023.
The second author was supported in part by UKRI grant MR/S032657/1.

\bibliographystyle{amsalpha}

\end{document}